\documentclass[12pt]{amsart}
\usepackage{amsfonts}
\usepackage{amssymb}
\numberwithin{equation}{section}
\theoremstyle{plain}
 \newtheorem{thm}{Theorem}[section]

 \newtheorem{prop}[thm]{Proposition}

\theoremstyle{definition}
 
 \newtheorem{ex}[thm]{Example}

\begin{document}
\setlength{\baselineskip}{18pt}
\setlength{\parindent}{1.8pc}
\title{Monotonicity and non-monotonicity of domains\\
of stochastic integral operators}
\author{Ken-iti Sato}
\dedicatory{To the memory of K. Urbanik}

\thanks{
\indent {\em Author's address:} 
Hachiman-yama 1101-5-103, Tenpaku-ku, Nagoya, 468-0074 Japan. 
{\em E-mail:} ken-iti.sato@nifty.ne.jp\;}
\maketitle

\begin{abstract}
A L\'evy process on $\mathbb{R}^d$ with distribution $\mu$ at time $1$ is denoted
by $X^{(\mu)}=\{X_t^{(\mu)}\}$. If the improper stochastic integral 
$\int_0^{\infty-} f(s)dX_s^{(\mu)}$
of $f$ with respect to $X^{(\mu)}$ is definable, its distribution is denoted
by $\Phi_f(\mu)$. The class of all 
infinitely divisible distributions
$\mu$ on $\mathbb{R}^d$ such that $\Phi_f(\mu)$ is definable is denoted by
$\mathfrak D(\Phi_f)$. 
The class $\mathfrak D(\Phi_f)$, its two extensions $\mathfrak{D}_{\mathrm{c}}(\Phi_f)$ 
and $\mathfrak{D}_{\mathrm{e}}(\Phi_f)$ (compensated and essential),
and its restriction $\mathfrak D^0(\Phi_f)$ (absolutely definable) are studied.
It is shown that $\mathfrak{D}_{\mathrm{e}}(\Phi_f)$
is monotonic with respect to $f$, which means that
$|f_2|\leqslant|f_1|$ implies $\mathfrak{D}_{\mathrm{e}}(\Phi_{f_1})\subset 
\mathfrak{D}_{\mathrm{e}}(\Phi_{f_2})$. Further,
$\mathfrak D^0(\Phi_f)$ is monotonic with respect to
$f$ but neither $\mathfrak D(\Phi_f)$ nor $\mathfrak{D}_{\mathrm{c}}(\Phi_f)$ is 
monotonic with respect to $f$. Furthermore, there exist $\mu$, $f_1$, and $f_2$ 
such that $0\leqslant f_2\leqslant f_1$,
$\mu\in\mathfrak D(\Phi_{f_1})$, and $\mu\not\in\mathfrak D(\Phi_{f_2})$.
An explicit example for this 
is related to some properties of a class of martingale L\'evy processes.
\end{abstract}

\section{Introduction and results}

Let $ID(\mathbb{R}^d)$ be the class of infinitely divisible distributions on the
$d$-dimensional Euclidean space $\mathbb{R}^d$.
For each $\mu\in ID(\mathbb{R}^d)$ let $X^{(\mu)}=\{ X_t^{(\mu)}, t\geqslant0\}$ be the
L\'evy process on $\mathbb{R}^d$ satisfying $\mathcal L(X_1^{(\mu)})=\mu$.  Here $\mathcal L(Y)$
denotes the distribution of $Y$ for any random element $Y$.  
Given $\mu\in ID(\mathbb{R}^d)$ and 
a real-valued measurable nonrandom function $f$ on $[0,\infty)$, we 
say, as in \cite{S05D}, that $f$ is {\it locally $X^{(\mu)}$-integrable} if the
stochastic integral $\int_B f(s)dX_s^{(\mu)}$ of $f$ with 
respect to $X^{(\mu)}$ is definable for each bounded Borel set $B$ in
$[0,\infty)$ in the sense of Urbanik and Woyczy\'nski \cite{UW67},\linebreak 
Rajput and Rosinski \cite{RR89}, Kwapie\'n and Woyczy\'nski \cite{KW92}, and
Sato \cite{S04,S05D,S05A}.  We write $\int_0^t f(s)dX_s^{(\mu)}=\int_{[0,t]}
 f(s)dX_s^{(\mu)}$. Since this is an additive process in law, we use an
additive process modification (see \cite{S} for terminology). 
For $\mu$ fixed, let
\begin{equation}\label{1}
\mathbf L(X^{(\mu)})=\{f\colon \text{$f$ is  locally $X^{(\mu)}$-integrable}\}.
\end{equation}
Characterization of $\mathbf L(X^{(\mu)})$ in terms of the L\'evy--Khintchine triplet
of $\mu$ is given in \cite{RR89,KW92,S04,S05D}.  It is known that 
$\mathbf L(X^{(\mu)})$
is a generalization of Orlicz spaces, 
one of whose properties is that $\mathbf L(X^{(\mu)})$
 is {\em monotonic}.  By this we 
mean that, if
\begin{equation}\label{3}
\text{$f_1$ and $f_2$ are measurable and }|f_2|\leqslant |f_1|,
\end{equation}
then $f_2\in\mathbf L(X^{(\mu)})$ whenever $f_1\in\mathbf L(X^{(\mu)})$.
Given $f$, denote
\begin{equation}\label{1a}
\mathbf D[f]=\mathbf D[f; \mathbb{R}^d]=\{\mu\in ID(\mathbb{R}^d)\colon 
\text{$f$ is  locally $X^{(\mu)}$-integrable}\}.
\end{equation}
Then \eqref{3} implies $\mathbf D[f_1]\subset\mathbf D[f_2]$.  We express this
property by saying that $\mathbf D[f]$ is {\it monotonic with respect to }$f$.

Let $\mu\in ID(\mathbb{R}^d)$.  We say that {\it the improper stochastic integral of $f$
with respect to
$X^{(\mu)}$ is definable} if  $f\in\mathbf L(X^{(\mu)})$ and if 
$\int_0^t f(s)dX_s^{(\mu)}$ is convergent in probability 
(equivalently, convergent almost surely) in $\mathbb{R}^d$ as $t\to\infty$. 
The limit is denoted by $\int_0^{\infty-} f(s)dX_s^{(\mu)}$.
This notation will help to distinguish it from the stochastic integral (with
random integrand in general) up to infinity of Cherny and Shiryaev \cite{CS05}.
We define
\begin{equation}\label{5a}
\Phi_f(\mu)=\mathcal L \left( \int_0^{\infty-} f(s)dX_s^{(\mu)}\right).
\end{equation}
Two extended notions
and one restricted notion of definability of improper stochastic integrals are
introduced in \cite{S05D,S05A}. 
We say that {\it the compensated improper stochastic integral of $f$ with respect to
$X^{(\mu)}$ is definable} if  $f\in\mathbf L(X^{(\mu)})$ and if there is $q\in\mathbb{R}^d$
 such that  $\int_0^{\infty-} f(s)dX_s^{(\mu*\delta_{-q})}$ is definable.
Here $\delta_{-q}$ is the distribution concentrated at $-q$. 
We say that {\it the essential improper integral of $f$ with respect to 
$X^{(\mu)}$ is definable} if  $f\in\mathbf L(X^{(\mu)})$ and if there is 
a nonrandom $\mathbb{R}^d$-valued
function $q_t$ on $[0,\infty)$ such that $\int_0^t f(s)dX_s^{(\mu)}-q_t$ is
convergent in probability in $\mathbb{R}^d$ as $t\to\infty$.
We say that {\it the improper integral 
$\int_0^{\infty-} f(s)dX_s^{(\mu)}$ is absolutely definable} if  
$f\in\mathbf L(X^{(\mu)})$ and if
\begin{equation}\label{7}
\int_0^{\infty} |C_{\mu}(f(s)z)| ds<\infty\text{ for all $z\in\mathbb{R}^d$}.
\end{equation}
Here $C_{\mu}(z)$ is the cumulant function of $\mu$, that is, the complex-valued
continuous function on $\mathbb{R}^d$ with $C_{\mu}(0)=0$ such that the characteristic
function $\widehat\mu(z)$ of $\mu$ is expressed as $\widehat\mu(z)=e^{C_{\mu}(z)}$.
For any measurable function $f$ on $[0,\infty)$, we denote
{\allowdisplaybreaks
\begin{align*}
\mathfrak{D}^0(\Phi_f)&=\mathfrak{D}^0(\Phi_f;\mathbb{R}^d)=\left\{ \mu\in ID(\mathbb{R}^d)
\colon \int_0^{\infty-} f(s)dX_s^{(\mu)}\text{ is absolutely definable}\right\},\\
\mathfrak D(\Phi_f)&=\mathfrak D(\Phi_f;\mathbb{R}^d)=\left\{ \mu\in ID(\mathbb{R}^d)
\colon \int_0^{\infty-} f(s)dX_s^{(\mu)}\text{ is definable}\right\},\\
\mathfrak{D}_{\mathrm{c}}(\Phi_f)&=\mathfrak{D}_{\mathrm{c}}(\Phi_f;\mathbb{R}^d)=\{\mu\in ID(\mathbb{R}^d)\colon \text{ compensated improper integral of 
$f$}\\
&\phantom{QQQQQQQQQQQQQQQQr}\text{ with respect to $X^{(\mu)}$ is definable}\},\\
\mathfrak{D}_{\mathrm{e}}(\Phi_f)&=\mathfrak{D}_{\mathrm{e}}(\Phi_f;\mathbb{R}^d)=\{\mu\in ID(\mathbb{R}^d)\colon \text{ essential improper integral of 
$f$}\\
&\phantom{QQQQQQQQQQQQQQQQr}\text{ with respect to $X^{(\mu)}$ is definable}\}.
\end{align*}
Further we denote, for} $\mu\in ID(\mathbb{R}^d)$, 
\[
\mathbf L^{\infty-}(X^{(\mu)})=\{f\colon f \text{ is measurable and }
\mu\in\mathfrak D(\Phi_f;\mathbb{R}^d)\}.
\]
It is known that
\begin{equation}\label{8}
\mathfrak{D}^0(\Phi_f)\subset \mathfrak D(\Phi_f)\subset \mathfrak{D}_{\mathrm{c}}(\Phi_f) \subset \mathfrak{D}_{\mathrm{e}}(\Phi_f).
\end{equation}

We are interested in the problem whether $\mathfrak D(\Phi_f)$, $\mathfrak{D}^0(\Phi_f)$, 
$\mathfrak{D}_{\mathrm{c}}(\Phi_f)$, and $\mathfrak{D}_{\mathrm{e}}(\Phi_f)$ are monotonic with respect to $f$. 
Clearly, $\mathfrak D(\Phi_f)$ is  monotonic with respect to $f$ if and only if,
for every $\mu\in ID(\mathbb{R}^d)$, $\mathbf L^{\infty-}(X^{(\mu)})$ is monotonic.  

Our results are the following.

\begin{thm}\label{t1}
The class $\mathfrak{D}_{\mathrm{e}}(\Phi_f)$ is monotonic with respect to $f$.
\end{thm}

\begin{thm}\label{t2}
The class $\mathfrak D^0(\Phi_f)$ is monotonic with respect to $f$.
\end{thm}

The class $\mathfrak D(\Phi_f)$ is not monotonic with respect to $f$. That is, 
for some $\mu\in ID(\mathbb{R}^d)$, $\mathbf L^{\infty-}(X^{(\mu)})$ is
 not monotonic. 
 In order to specify $\mu$, we
 use the L\'evy--Khintchine triplet $(A,\nu,\gamma)$
of $\mu\in ID(\mathbb{R}^d)$ in the sense that
\begin{equation*}
C_{\mu}(z)= -\frac12 \langle z,Az\rangle +\int_{\mathbb{R}^d}\left(
e^{i\langle z,x\rangle}-1-\frac{i\langle z,x\rangle}{1+|x|^2}\right)\nu(dx)
+i\langle\gamma,z\rangle,
\end{equation*}
where $A$ is a $d\times d$ symmetric nonnegative-definite matrix, 
called the Gaussian covariance matrix of $\mu$, $\nu$ is a
measure on $\mathbb{R}^d$ satisfying $\nu(\{0\})=0$ and $\int_{\mathbb{R}^d}(|x|^2\land1)\nu(dx)
<\infty$, called the L\'evy measure of $\mu$, and $\gamma$ is an element 
of $\mathbb{R}^d$, called the location parameter of $\mu$. Sometimes we denote 
$\mu=\mu_{(A,\nu,\gamma)}$. We say that a measure $\rho$ on $\mathbb{R}^d$
is {\it symmetric} if $\rho(B)=\rho(-B)$
for all Borel sets $B$.

\begin{thm}\label{t3}
Let $\mu=\mu_{(A,\nu,\gamma)}
\in ID(\mathbb{R}^d)$ with $A$ arbitrary and $\nu$ symmetric.

{\rm(i)}  If $f_1$ and $f_2$ satisfy \eqref{3} and if $\mu\in\mathfrak{D}_{\mathrm{c}}(\Phi_{f_1})$,
then $\mu\in\mathfrak{D}_{\mathrm{c}}(\Phi_{f_2})$.

{\rm(ii)}  Assume that $\gamma=0$.
If $f_1$ and $f_2$ satisfy \eqref{3} and if $\mu\in\mathfrak D(\Phi_{f_1})$,
then $\mu\in\mathfrak D(\Phi_{f_2})$. That is, $\mathbf L^{\infty-}(X^{(\mu)})$ 
is monotonic.

{\rm(iii)}  Assume that $\gamma\ne0$ and $\int_{|x|>1} |x|\nu(dx)<\infty$.
Then $\mathbf L^{\infty-}(X^{(\mu)})$ is not monotonic.
\end{thm}

A simple example for Theorem \ref{t3} (iii)
 is the case where $X^{(\mu)}$ is a Brownian 
motion with drift.

We ask a question whether there exist $\mu$, $f_1$, and $f_2$ such that 
$0\leqslant f_2\leqslant f_1$, $f_1\in\mathbf L^{\infty-}(X^{(\mu)})$, and 
$f_2\not\in\mathbf L^{\infty-}(X^{(\mu)})$.
The next theorem gives more than the affirmative answer.

\begin{thm}\label{t4a} Let $f_1(s)$ be a real-valued function which vanishes
on $[0,a)$ and is continuous on $[a,\infty)$ with some $a\geqslant0$.  Let
$\mu\in\mathfrak D(\Phi_{f_1})\setminus\mathfrak D^0(\Phi_{f_1})$.  Then there is a
nonrandom open set $D$ in $[a,\infty)$ such that $\mu\not\in\mathfrak D(\Phi_{f_2})$
for $f_2(s)=f_1(s)1_D(s)$.
\end{thm}

Notice that this theorem and Theorem \ref{t2} give a characterization of the
property that $\mathfrak D(\Phi_{f_1})\setminus\mathfrak D^0(\Phi_{f_1})\neq\emptyset$.

We say that $f(s)\asymp g(s)$ as $s\to\infty$ if there are positive constants
$c_1$ and $c_2$ such that $0<c_1 f(s)\leqslant g(s)\leqslant c_2 f(s)$ for all large
$s$.

\begin{ex}\label{e1}
Let $f_1(s)$ be a locally square-integrable function on $[0,\infty)$. Suppose
that $f_1(s)\asymp s^{-1}$ as $s\to\infty$ and that
there are positive constants $c$ and $s_0$ such that
\(
\int_{s_0}^{\infty} |f_1(s)-c s^{-1}|ds<\infty.
\)
Then Theorem 2.8 of \cite{S05A} says that the class  
$\mathfrak D(\Phi_{f_1})\setminus\mathfrak D^0(\Phi_{f_1})$ is nonempty and that
$\mu=\mu_{(A,\nu,\gamma)}\in \mathfrak D(\Phi_{f_1})\setminus\mathfrak D^0(\Phi_{f_1})$ 
if and only if $\int_{|x|>1} |x|\nu(dx)<\infty$, $\int_{\mathbb{R}^d} x\mu(dx)=0$,
${\displaystyle\lim_{t\to\infty}} \int_{s_0}^t s^{-1}ds\int_{|x|>s} x\nu(dx)$ exists
in $\mathbb{R}^d$, and $\int_{s_0}^{\infty}s^{-1}\left|
\int_{|x|>s}x\nu(dx)\right|ds=\infty$.   Distributions satisfying these
conditions will be given in Example \ref{e2}.
\end{ex}

We show that the class $\mathfrak{D}_{\mathrm{c}}(\Phi_f)$ is not monotonic with respect to
$f$.

\begin{thm}\label{t5a}
Let  $f_1(s)=s^{-1}1_{[1,\infty)}(s)$. 
Suppose that $\mu\in\mathfrak D(\Phi_{f_1})$ and that the L\'evy measure $\nu$
of $\mu$ satisfies
\begin{equation}\label{t5a.1}
\left|\int_{|x|>s}x_j\nu(dx)\right|\sim\frac{c}{\log s}\quad\text{as }
s\to\infty
\end{equation}
for some $j\in\{1,\ldots,d\}$ and $c>0$. Then there is a 
nonrandom open set $D$ in $[1,\infty)$ such that
$\mu\not\in\mathfrak{D}_{\mathrm{c}}(\Phi_{f_2})$  for $f_2(s)=f_1(s)1_D(s)$.
\end{thm}

Here $x_j$ is the $j$th coordinate of $x\in\mathbb{R}^d$. In Theorem \ref{t5a}
recall that $\mu\in\mathfrak D(\Phi_{f_1})$ implies $\int_{|x|>1}|x|\nu(dx)<\infty$ 
by virtue of Theorem 2.8 of \cite{S05A}.

\begin{ex}\label{e2}
In Example 2.9 of \cite{S05A} we have introduced
 the measure $\nu$ concentrated on 
$\{x\in\mathbb{R}^d\colon |x|=2,3,\ldots\}$ given by
\begin{equation*}
\nu(B)=\int_{S_0}\lambda(d\xi)\sum_{n\in\mathbb{Z}} 
1_B (n\xi)a_n\quad\text{for Borel sets $B$},
\end{equation*}
where $S_0$ is a nonempty Borel set on the unit sphere 
$\{|\xi|=1\}$ satisfying
$S_0\cap(-S_0)=\emptyset$, $\lambda$ is a finite measure 
on $S_0$ satisfying
$\int_{S_0} \xi\lambda(d\xi)\ne 0$, $\mathbb{Z}$ is the 
class of all integers, and 
$a_n$, $n\in\mathbb{Z}$, are such that $a_0=a_1=a_{-1}=0$ and,
for positive integers $n$, $m$,
{\allowdisplaybreaks
\begin{align*}
a_n&=\frac1{n}\left(\frac1{\log n}-\frac1{\log(n+1)}
\right),\;a_{-n}=0
 \quad\text{if $2^{m^2}<n<2^{(m+1)^2}$, $m$ odd,}\\
a_n&=0,\;a_{-n}=\frac1{n}\left(\frac1{\log n}-\frac1{\log(n+1)}
\right)
\quad\text{if $2^{m^2}<n<2^{(m+1)^2}$, $m$ even,}\\
a_n&=\frac1{n}\left(\frac1{\log n}+
\frac1{\log(n+1)}\right),\;a_{-n}=0
 \quad\text{if $n=2^{m^2}$, $m$ even,}\\
a_n&=0,\;a_{-n}=\frac1{n}\left(\frac1{\log n}+
\frac1{\log(n+1)}\right)
\quad\text{if $n=2^{m^2}$, $m$ odd.}
\end{align*}
It is shown that} $\sum_{|n|\geqslant2} |n| a_n<\infty$
and that, for $k=3,4,\ldots$,
\begin{equation*}
\sum_{|n|\geqslant k}n a_n=\begin{cases}
(\log k)^{-1} \quad & \text{ if $2^{m^2}<k\leqslant 2^{(m+1)^2}$, $m$ odd,}\\
-(\log k)^{-1} \quad &\text{ if $2^{m^2}<k\leqslant 2^{(m+1)^2}$, $m$ even.}
\end{cases}
\end{equation*}
Thus
{\allowdisplaybreaks
\begin{gather*}
\int_{|x|>1}|x|\nu(dx)=\int_{S_0}\lambda(d\xi)
\sum_{|n|\geqslant2}|n|a_n<\infty,\\
\int_1^{\infty}s^{-1}ds\left|\int_{|x|>s}x\nu(dx)\right|=\left|
\int_{S_0}\xi\lambda(d\xi)\right|\,\int_1^{\infty}s^{-1}ds\left|
\sum_{|n|>s}
na_n\right|=\infty,\\
\int_{|x|>s}x_j\nu(dx)=\int_{S_0}\xi_j\lambda(d\xi)\sum_{|n|>s}na_n.
\end{gather*}
Further it is shown that} 
$\int_1^t s^{-1}ds\int_{|x|>s}x\nu(dx)$ 
is convergent as $t\to\infty$. 
Let $\mu=\mu_{(A,\nu,\gamma)}$ with $\gamma=-\int_{\mathbb{R}^d}x|x|^2(1+|x|^2)^{-1}
\nu(dx)$ and $A$ arbitrary.  Then $\int_{\mathbb{R}^d}x\mu(dx)=0$ and
$\mu\in\mathfrak D(\Phi_{f_1})\setminus\mathfrak D^0(\Phi_{f_1})$ for
$f_1(s)=s^{-1}1_{[1,\infty)}(s)$, since the conditions stated in Example
\ref{e1} are satisfied.
Choosing $j\in\{1,\ldots,d\}$ such that $\int_{S_0}\xi_j\lambda(d\xi)\neq0$,
we can apply Theorem \ref{t5a} to this distribution $\mu$.
We can also apply Theorem \ref{t4a} to this $f_1$ and this $\mu$.
If $A=0$, then the process 
$X^{(\mu)}$ is a compensated compound Poisson process.  

Using the measure $\nu$ above, consider
\[
\widetilde\nu(B)=\nu(B)+\frac{1}{2\log 2}\int_{S_0}\lambda(d\xi)1_B(2\xi)
\quad\text{for Borel sets $B$}
\]
and define $\widetilde\mu\in ID(\mathbb{R}^d)$ by
\[
C_{\widetilde\mu}(z)=\int_{\mathbb{R}^d}(e^{i\langle z,x\rangle}-1)\widetilde\nu(dx).
\]
Then 
\[
\int_{\mathbb{R}^d}x\widetilde\mu(dx)=\int_{\mathbb{R}^d}x\widetilde\nu(dx)=\int_{\mathbb{R}^d}x\nu(dx)+
\frac{1}{\log 2}\int_{S_0}\xi\lambda(d\xi)=0.
\]
The distribution $\widetilde\mu$ also belongs to
$\mathfrak D(\Phi_{f_1})\setminus\mathfrak D^0(\Phi_{f_1})$ for
$f_1(s)=s^{-1}1_{[1,\infty)}(s)$ and Theorem \ref{t5a} applies to
$\widetilde\mu$ by the same reason as for $\mu$. Theorem \ref{t4a} also
applies to $f_1(s)=s^{-1}1_{[1,\infty)}(s)$ and $\widetilde\mu$.
The L\'evy process $X^{(\widetilde\mu)}$ associated with $\widetilde\mu$ is a 
compound Poisson process with mean $0$.  
\end{ex}

In Section 2 we will give proofs of all theorems stated above. 
The process $X^{(\mu)}$ associated with $\mu$ in Theorem \ref{t5a} 
is a martingale L\'evy process and the processes
$\int_0^t f_1(s)dX_s^{(\mu)}$ and $\int_0^t f_2(s)dX_s^{(\mu)}$ 
have intriguing properties, which we will discuss in Section 3. 
Applications of Theorems \ref{t1} and \ref{t2} 
to some types of $f$ will be given in Section 4.
Determination of $\mathfrak D(\Phi_f)$ for some $f$ is made.

\section{Proofs}

In the following three propositions let $\mu=\mu_{(A,\nu,\gamma)}\in ID(\mathbb{R}^d)$
and $f\in \mathbf L(X^{(\mu)})$.  We present necessary and sufficient conditions for
$\mu$ to belong to $\mathfrak D(\Phi_f)$, $\mathfrak D^0(\Phi_f)$, or $\mathfrak{D}_{\mathrm{e}} (\Phi_f)$.

\begin{prop}\label{p1}
The following three statements are equivalent.

{\rm(a)}\quad $\mu\in\mathfrak D(\Phi_f)$.

{\rm(b)}\quad $\int_0^t C_{\mu}(f(s)z)ds$ is convergent in $\mathbb{C}$ as $t\to\infty$
for each $z\in\mathbb{C}$.

{\rm(c)}\quad $\mu$ satisfies the following:
{\allowdisplaybreaks
\begin{gather}
\int_0^{\infty} f(s)^2 (\mathrm{tr}\, A)ds<\infty,\label{2.1}\\
\int_0^{\infty} ds \int_{\mathbb{R}^d}(|f(s)x|^2\land1)\nu(dx)<\infty,\label{2.2}\\
\begin{split}
&\int_0^t f(s)\left(\gamma+\int_{\mathbb{R}^d}x\left(\frac{1}{1+|f(s)x|^2}-\frac{1}{1+|x|^2}
\right)\nu(dx) \right) ds\\
&\text{ is convergent in $\mathbb{R}^d$ as $t\to\infty$.}
\end{split}\label{2.3}
\end{gather}}
\end{prop}
\indent {\it Proof}.  See Proposition 5.5 of \cite{S05D} and Propositions 2.2 and 2.6
of \cite{S05A}. It follows from $f\in\mathbf L(X^{(\mu)})$ that
$\int_0^t |C_{\mu}(f(s)z)|ds<\infty$ and that
\[
\int_0^t \left|f(s)\left(\gamma+\int_{\mathbb{R}^d}x\left(\frac{1}{1+|f(s)x|^2}-\frac{1}{1+|x|^2}
\right)\nu(dx) \right)\right| ds<\infty
\]
for $t\in(0,\infty)$, as is shown in Proposition 2.17 and Corollary 2.19 of
\cite{S05D}. \qed

\begin{prop}\label{p2}
A distribution $\mu$ is in $\mathfrak{D}_{\mathrm{e}} (\Phi_f)$ if and only if \eqref{2.1} and \eqref{2.2} 
are satisfied.
\end{prop}

{\it Proof}. See  Proposition 5.6 of \cite{S05D} or Proposition 2.6 of \cite{S05A}.
\qed

\begin{prop}\label{p3}
A distribution $\mu$ is in $\mathfrak D^0 (\Phi_f)$ if and only if \eqref{2.1}, \eqref{2.2},
and
\begin{equation}\label{2.4}
\int_0^{\infty}
 \left|f(s)\left(\gamma+\int_{\mathbb{R}^d}x\left(\frac{1}{1+|f(s)x|^2}-\frac{1}{1+|x|^2}
\right)\nu(dx) \right)\right| ds<\infty.
\end{equation}
\end{prop}

{\it Proof}.  For fixed $u\in\mathbb{R}$ denote by $\mu^u$ a probability measure such that
$\mu^u(B)=\int 1_B(ux)\mu(dx)$ for all Borel sets $B$. Let $(A^u,\nu^u,
\gamma^u)$ be the triplet of $\mu^u$. Then $A^u=u^2 A$, $\nu^u(B)=\int 1_B(ux)\nu(dx)$, 
and
\[
\gamma^u=u\gamma+\int_{\mathbb{R}^d}ux\left(\frac{1}{1+|ux|^2}-\frac{1}{1+|x|^2}
\right)\nu(dx).
\]
Notice that
\begin{equation}\label{2.4a}
\int_{\mathbb{R}^d}|ux|\left|\frac{1}{1+|ux|^2}-\frac{1}{1+|x|^2}
\right|\nu(dx)\leqslant\int_{\mathbb{R}^d}\frac{|ux|(|x|^2+|ux|^2)}{(1+|ux|^2)(1+|x|^2)}
\nu(dx)<\infty.
\end{equation}
Let
\begin{equation}\label{2.5}
\varphi(u)=\mathrm{tr}\, A^u +\int_{\mathbb{R}^d}(|x|^2\land1)\nu^u(dx)+|\gamma^u|.
\end{equation}
When $u=f(s)$, $\mu^u$ and $(A^u,\nu^u,
\gamma^u)$ are written as $\mu^{f(s)}$ and $(A^{f(s)},\nu^{f(s)},
\gamma^{f(s)})$. 
The properties \eqref{2.1}, \eqref{2.2},
and \eqref{2.4} combined are expressed by
\begin{equation}\label{2.6}
\int_0^{\infty} \varphi(f(s))ds<\infty.
\end{equation}
We note that 
\begin{align*}
|C_{\mu}( f(s)z)|&=|C_{\mu^{f(s)}}(z)|\\
&\leqslant\frac{|z|^2}{2} \mathrm{tr}\, A^{f(s)}+3(1+|z|^2)
\int_{\mathbb{R}^d}(|x|^2\land1)\nu^{f(s)}(dx)+|z||\gamma^{f(s)}|
\end{align*}
(see, in \cite{S05D}, (2.5)--(2.7) and line 3 of the proof of Theorem 3.14).
 Hence, if \eqref{2.1}, \eqref{2.2},
and \eqref{2.4} are satisfied, then \eqref{7} is satisfied, that is, 
$\mu\in\mathfrak D^0(\Phi_f)$. 

Conversely, assume that $\mu\in\mathfrak D^0(\Phi_f)$. Then $\mu\in\mathfrak D(\Phi_f)$ and
\eqref{2.1} and \eqref{2.2} follow from Proposition \ref{p1}.  We have
\[
\mathrm{Im}\,C_{\mu}(f(s)z)=\mathrm{Im}\,C_{\mu^{f(s)}}(z)=\int_{\mathbb{R}^d}\left(\sin\langle z,x\rangle-
\frac{\langle z,x\rangle}{1+|x|^2}\right)\nu^{f(s)}(dx) +\langle \gamma^{f(s)},z\rangle.
\]
For fixed $z$,
\[
\sin\langle z,x\rangle-\frac{\langle z,x\rangle}{1+|x|^2}=\begin{cases}
O(|x|^3),\quad & |x|\to0,\\
O(1),\quad & |x|\to\infty.
\end{cases}
\]
Hence it follows from \eqref{2.2} that
\begin{align*}
\int_{\mathbb{R}^d}\left|\sin\langle z,x\rangle-
\frac{\langle z,x\rangle}{1+|x|^2}\right|\nu^{f(s)}(dx)&\leqslant c_z 
\int_{\mathbb{R}^d}(|x|^2\land1)\nu^{f(s)}(dx)\\
&=c_z\int_{\mathbb{R}^d}(|f(s)x|^2\land1)\nu(dx),
\end{align*}
where $c_z$ is a constant depending on $z$. Thus we obtain
\[
\int_0^{\infty}|\langle \gamma^{f(s)},z\rangle| ds<\infty
\]
from $\int_0^{\infty} |\mathrm{Im}\,C_{\mu}(f(s)z)| ds<\infty$.
Choosing $z=(\delta_{jk})_{1\leqslant k\leqslant d}$, $1\leqslant j\leqslant d$, we obtain \eqref{2.4}.
\qed

\medskip
{\it Proof of Theorem \ref{t1}}. Use Proposition \ref{p2}. 
Let $f_1$ and $f_2$ satisfy \eqref{3}. 
Suppose that $\mu\in\mathfrak{D}_{\mathrm{e}}(\Phi_{f_1})$.  Then \eqref{2.1} and \eqref{2.2} hold with 
$f_1$ in place of $f$.  Since $|f_2|\leqslant|f_1|$, it follows that 
\eqref{2.1} and \eqref{2.2} hold with 
$f_2$ in place of $f$.  This means that $\mu\in\mathfrak{D}_{\mathrm{e}}(\Phi_{f_2})$.  \qed

\medskip
{\it Proof of Theorem \ref{t2}}. Use Proposition \ref{p3}. 
Let $f_1$ and $f_2$ satisfy \eqref{3} and suppose that 
$\mu\in\mathfrak D^0(\Phi_{f_1})$.  
Using the function $\varphi(u)$ in \eqref{2.5} induced by $\mu=\mu_{(A,\nu,\gamma)}$,
we have
\begin{equation}\label{2.7}
\int_0^{\infty}\varphi(f_1(s))ds<\infty.
\end{equation}
Let us use
\begin{equation}\label{2.8}
\widetilde\varphi(u)=\mathrm{tr}\, A^u +\int_{\mathbb{R}^d}(|x|^2\land1)\nu^u(dx)+
\sup_{v\in\mathbb{R},\,|v|\leqslant|u|}|\gamma^v|.
\end{equation}
We have 
$\varphi(u)\leqslant\widetilde\varphi(u)\leqslant(3/2) \varphi(u)$ as in  Proposition 3.10 of \cite{S05D}.
Thus \eqref{2.7} is equivalent to
\begin{equation}\label{2.9}
\int_0^{\infty} \widetilde\varphi(f_1(s))ds<\infty.
\end{equation}
The function $\widetilde\varphi$ enjoys the property
that $\widetilde\varphi(f_2(s))\leqslant\widetilde\varphi(f_1(s))$ whenever $|f_2(s)|\leqslant |f_1(s)|$. Hence, 
$\int_0^{\infty}\widetilde\varphi(f_2(s))ds<\infty$. This means $\mu\in\mathfrak D^0(\Phi_{f_2})$.  
\qed

\medskip
{\it Proof of Theorem \ref{t3}}.  (i)  Let $f_1$ and $f_2$ satisfy \eqref{3}.
Assume that $\mu\in\mathfrak{D}_{\mathrm{c}}(\Phi_{f_1})$.  Then there is $q\in\mathbb{R}^d$ such that 
$\mu*\delta_{-q}\in\mathfrak D(\Phi_{f_1})$.  Thus $\int_0^{\infty} f_1(s)^2 A ds<\infty$,
$\int_0^{\infty}ds\int_{\mathbb{R}^d}(|f_1(s)x|^2\land1) \nu(dx)<\infty$, and
$\int_0^t f_1(s)(\gamma-q)ds$ is convergent, since $\nu$ is symmetric.  We may and do
choose $q=\gamma$. Then we see that 
$\mu*\delta_{-q}\in\mathfrak D(\Phi_{f_2})$. It follows that
$\mu\in\mathfrak{D}_{\mathrm{c}}(\Phi_{f_2})$.

(ii)  Look back to the argument above with $\gamma=0$. Then the proof is
evident.

(iii) Let
\[
f_1(s)=\begin{cases}0 \quad & \text{if $0\leqslant s<1$}\\
s^{-1} \quad & \text{if $n\leqslant s<n+1$ with $n$ odd}\\
-s^{-1} \quad & \text{if $n\leqslant s<n+1$ with $n$ even} \end{cases}
\]
and let $f_2(s)=s^{-1}1_{[1,\infty)}(s)$.  Then $|f_2|=|f_1|$.  Applying 
Theorem 2.8 of \cite{S05A}, we see that $f_2\not\in \mathbf L^{\infty-}(X^{(\mu)})$
since $\gamma\neq0=\int_{\mathbb{R}^d}x|x|^2(1+|x|^2)^{-1}\nu(dx)$. 
On the other hand, $f_1\in \mathbf L^{\infty-}(X^{(\mu)})$ 
by virtue of Proposition \ref{p1}. Indeed,
$\int_0^{\infty}ds\int_{\mathbb{R}^d}(|f_1(s)x|^2\land1)\nu(dx)<\infty$ by the same
reasoning as in the proof of Lemma 2.7 of \cite{S05A}, and
\[
\int_0^t f_1(s)\left(\gamma+\int_{\mathbb{R}^d}x\left(\frac{1}{1+|f_1(s)x|^2}-\frac{1}{1+|x|^2}
\right)\nu(dx) \right) ds=\int_0^t f_1(s)ds\gamma,
\]
which is convergent in $\mathbb{R}^d$ as $t\to\infty$. Hence $\mathbf L^{\infty-}(X^{(\mu)})$
is not monotonic.  \qed

\medskip
{\it Proof of Theorem \ref{t4a}}. 
Let $(A,\nu,\gamma)$ be the triplet of $\mu$. We use an $\mathbb{R}^d$-valued
function
\begin{equation}\label{2.9a}
h(s)=f_1(s)\gamma+\int_{\mathbb{R}^d} f_1(s)x\left(\frac{1}{1+|f_1(s)x|^2}
-\frac{1}{1+|x|^2}\right)\nu(dx).
\end{equation}
Using \eqref{2.4a},
we see that $h(s)$ is continuous on $[a,\infty)$. 
Since $\mu\in\mathfrak D(\Phi_{f_1})$, we have $\int_a^{\infty} f_1(s)^2 (\mathrm{tr}\, A)
 ds<\infty$,
$\int_a^{\infty} ds\int_{\mathbb{R}^d}(|f_1(s)x|^2\land1)\nu(dx)<\infty$, and
$\int_a^t h(s)ds$ is convergent in $\mathbb{R}^d$ as $t\to\infty$ (Proposition \ref{p1}).
Since $\mu\not\in\mathfrak D^0(\Phi_{f_1})$, we have $\int_a^{\infty} | h(s)|ds=\infty$
(Proposition \ref{p3}).
Choose and fix $j\in\{1,\ldots,d\}$ such that $\int_a^{\infty}|h_j(s)|ds=\infty$,
where $h_j(s)$ is the $j$th coordinate of $h(s)$.
Define $D^+=\{s\geqslant a\colon h_j(s)>0\}$, $D^-=\{s\geqslant a\colon h_j(s)<0\}$, 
and $D^0=\{s\geqslant a\colon h_j(s)=0\}$.  Then $D^+$ and $D^-$ are open in $[a,\infty)$.
Let $h_j^+(s)=h_j(s)\lor0$ and $h_j^- (s)=h_j^+ (s)-h_j(s)$.  
We see
that $\int_a^{\infty} h_j^+ (s)ds=\infty$ and $\int_a^{\infty} h_j^- (s)ds=\infty$. 
Let $D=D^+$ or $D^-$ (either will do). Let $f_2(s)=f_1(s) 1_D(s)$.  
Then
{\allowdisplaybreaks
\begin{align*}
&\int_0^t \left(f_2(s)\gamma +\int_{\mathbb{R}^d}f_2(s)x \left(
\frac{1}{1+|f_2(s)x|^2}-\frac{1}{1+|x|^2}\right) \nu(dx)\right)ds\\
&\quad=\int_a^t f_1(s)1_D (s)\left(\gamma +\int_{\mathbb{R}^d} x \left(
\frac{1}{1+|f_1(s)x|^2}-\frac{1}{1+|x|^2}\right) \nu(dx)\right)ds\\
&\quad=\int_a^t 1_D (s)h(s)ds.
\end{align*}
If $D=D^+$, then $\int_a^t 1_D(s)h_j(s)ds=\int_a^t h_j^+ (s)ds
\to\infty$ as $t\to\infty$.
If $D=D^-$, then $\int_a^t 1_D(s)h_j(s)ds=-\int_a^t h_j^- (s)ds
\to-\infty$ as $t\to\infty$.
Hence $\int_a^t 1_D (s)h(s)ds$ is not convergent in $\mathbb{R}^d$.
Hence $\mu\not\in\mathfrak D(\Phi_{f_2})$ by virtue of Proposition \ref{p1}.
\qed

\medskip
{\it Proof of Theorem \ref{t5a}}. Let the triplet of $\mu$ be $(A,\nu,\gamma)$.
In order to prove that $\mu\not\in\mathfrak{D}_{\mathrm{c}}(\Phi_{f_2})$ it is enough to show
that, for every $q\in\mathbb{R}^d$,
\[
\int_1^ts^{-1}1_D(s)\left( \gamma-q+\int_{\mathbb{R}^d}x\left(\frac{1}{1+|s^{-1}x|^2}
-\frac{1}{1+|x|^2}\right)\nu(dx)\right)ds
\]
is not convergent in $\mathbb{R}^d$ as $t\to\infty$.  Let $h(s)$ be as in \eqref{2.9a}.
This is an $\mathbb{R}^d$-valued function, continuous on $[1,\infty)$.
Since $\gamma=-\int_{\mathbb{R}^d}x|x|^2(1+|x|^2)^{-1}\nu(dx)$ (see Theorem
2.8 of \cite{S05A}), we have
\[
h(s)=\int_{\mathbb{R}^d} s^{-1}x\left(\frac{1}{1+|s^{-1}x|^2}-1\right)\nu(dx)
\quad\text{for }s\geqslant1.
\]
Choose $j$ as in \eqref{t5a.1} and let $h_j(s)$ be the $j$th coordinate
of $h(s)$.  Since $\int_1^t h(s) ds$ is convergent in $\mathbb{R}^d$ as $t\to\infty$ 
(see Theorem 2.8 of \cite{S05A}), $\int_1^t h_j(s) ds$ is convergent in $\mathbb{R}$.
We claim that
\begin{equation}\label{t5a.2}
\int_1^t |h_j(s)|ds\sim c\log\log t,\qquad t\to\infty.
\end{equation}
Indeed,
\[
\int_1^t |h_j(s)|ds\leqslant \int_1^ts^{-1}\left|\int_{|x|>s} x_j \nu(dx)\right| ds 
+I_1+I_2,
\]
where
\[
I_1=\int_1^ts^{-1}\left| \int_{|x|>s} \frac{x_j\nu(dx)}{1+|s^{-1}x|^2}\right|ds,
\quad I_2=\int_1^ts^{-1}\left| \int_{1<|x|\leqslant s} 
\frac{x_j|s^{-1}x|^2\nu(dx)}{1+|s^{-1}x|^2}\right|ds.
\]
We will denote positive constants by $c_1, c_2, \ldots$.
The quantities $I_1$ and $I_2$ are bounded in $t$, since
{\allowdisplaybreaks
\begin{align*}
I_1&\leqslant\int_1^{\infty} s^{-1}ds\int_{|x|>s}\frac{|x|\nu(dx)}{1+|s^{-1}x|^2}
=\int_{|x|>1}|x|\nu(dx)\int_1^{|x|}\frac{s^{-1}ds}{1+|s^{-1}x|^2}\\
&=\int_{|x|>1}|x|\left( \log|x|-\frac12\log(1+|x|^2)+\frac12\log 2\right)
\nu(dx)\\
&\leqslant c_1\int_{|x|>1} |x|\nu(dx)
\end{align*}
and since
\begin{align*}
I_2&\leqslant\int_1^{\infty} s^{-1}ds\int_{1<|x|\leqslant s}
\frac{|x|\,|s^{-1}x|^2\nu(dx)}{1+|s^{-1}x|^2}=\int_{|x|>1}|x|\nu(dx)
\int_{|x|}^{\infty}\frac{s^{-1}|s^{-1}x|^2ds}{1+|s^{-1}x|^2}\\
&=\frac{\log2}{2} \int_{|x|>1} |x|\nu(dx).
\end{align*}
Thus we obtain, for some} $s_0>1$,
\[
\int_1^t |h_j(s)|ds\leqslant c_2\int_{s_0}^t\frac{ds}{s\log s}+c_3
\leqslant c_2 \log\log t +c_4
\]
from condition \eqref{t5a.1}. Similarly,
\[
\int_1^t |h_j(s)|ds\geqslant c_5\int_{s_0}^t\frac{ds}{s\log s}-c_6
\geqslant c_5 \log\log t -c_7.
\]
Looking back more carefully, we see that \eqref{t5a.2} holds. We have,
a fortiori, $\int_1^{\infty}|h_j(s)|ds=\infty$.

Now define $D^+$, $D^-$, $D^0$, $h_j^+(s)$, and $h_j^-(s)$
 as in the proof of Theorem \ref{t4a}. Then $\int_1^{\infty}h_j^+(s)ds
 =\infty$ and $\int_1^{\infty}h_j^-(s)ds =\infty$.
We have
\begin{equation}\label{t5a.3}
\int_2^{\infty} 1_{D^0}(s)\frac{ds}{s\log s}<\infty,
\end{equation}
because it follows from
\[
|h_j(s)|\geqslant s^{-1}\left|\int_{|x|>s} x_j\nu(dx)\right|-\int_{|x|>s}
\frac{s^{-1}|x|\nu(dx)}{1+|s^{-1}x|^2} -\int_{|x|>s}
\frac{|s^{-1}x|^3\nu(dx)}{1+|s^{-1}x|^2} 
\]
that
{\allowdisplaybreaks
\begin{align*}
0&=\int_0^{\infty}1_{D^0}(s)|h_j(s)|ds\geqslant \int_1^{\infty}1_{D^0}(s)
\left|\int_{|x|>s} x_j\nu(dx)\right|\frac{ds}{s}-c_8\\
&\geqslant c_9\int_{s_0}^{\infty}1_{D^0}(s)\frac{ds}{s\log s}-c_{10}.
\end{align*}
Using \eqref{t5a.3},} we see that
\[
\limsup_{t\to\infty}\frac{1}{\log\log t}\int_2^t 1_{D^+}(s)
\frac{ds}{s\log s}+
\limsup_{t\to\infty}\frac{1}{\log\log t}\int_2^t 1_{D^-}(s)
\frac{ds}{s\log s}\geqslant 1.
\]
Choose $D=D^+$ or $D^-$ in such a way that
\[
\limsup_{t\to\infty}\frac{1}{\log\log t}\int_2^t 1_{D}(s)
\frac{ds}{s\log s}>0.
\]
Choose $t_n\to\infty$ such that $(\log\log t_n)^{-1}\int_2^{t_n}
1_D(s)(s\log s)^{-1}ds$ tends to some $b>0$. Then
\begin{equation}\label{t5a.4}
\frac{1}{\log\log t_n}\int_2^{t_n}1_D(s)\frac{ds}{s}\to\infty,
\end{equation}
since, for any $k>0$,
\[
\frac{1}{\log\log t_n}\int_2^{t_n}1_D(s)\frac{ds}{s}\geqslant
\frac{\log k}{\log\log t_n}\int_k^{t_n}1_D(s)\frac{ds}{s\log s}
\to b\log k.
\]
We claim that, for any choice of $q\in\mathbb{R}^d$, $\int_2^{t_n}1_D(s)
(h_j(s)-s^{-1}q_j)ds$ is divergent as $n\to\infty$.  
If $q_j=0$, then this is divergent since
\[
\int_2^{t_n} 1_D(s)h_j(s)ds=\int_2^{t_n}h_j^+(s)ds\quad\text{or}
\quad \int_2^{t_n}h_j^-(s)ds,
\]
which diverges to $\infty$ or to $-\infty$.
If $q_j\not=0$, then
\[
\int_2^{t_n}1_D(s)h_j(s)ds-q_j\int_2^{t_n}1_D(s)s^{-1}ds
\]
is divergent, because 
\[
\int_2^{t_n}1_D(s)|h_j(s)|ds\leqslant \int_2^{t_n}|h_j(s)|ds\sim c
\log\log t_n
\]
from \eqref{t5a.2} and because of \eqref{t5a.4}. The
proof is complete.  \qed

\section{Remarks on martingale L\'evy processes}

We have the following general result. Recall that
$f\in\mathbf L(X^{(\mu)})$ for all $\mu\in ID(\mathbb{R}^d)$ if and only if
$f$ is locally square-integrable on $[0,\infty)$ (see \cite{S05D}).

\begin{prop}\label{mart}
Let $X^{(\mu)}$ be a martingale L\'evy process on $\mathbb{R}^d$. Let $f$ be
locally square-integrable on $[0,\infty)$.  Then $Y_t=\int_0^t f(s)
dX_s^{(\mu)}$ is a martingale additive process.
\end{prop}

{\it Proof}. Let $\mu=\mu_{(A,\nu,\gamma)}$. 
We have $E|X_t^{(\mu)}|<\infty$ and
\[
0=EX_t^{(\mu)}=t\left(\gamma+\int_{\mathbb{R}^d}\frac{x|x|^2}{1+|x|^2}\nu(dx)\right).
\]
Let $(A_t,\nu_t,\gamma_t)$ be the triplet of $Y_t$. Proposition 2.6 of 
\cite{S05A} says that
{\allowdisplaybreaks
\begin{gather*}
A_t=\int_0^t f(s)^2 ds\,A,\\
\nu_t(B)=\int_0^t ds \int_{\mathbb{R}^d} 1_B(f(s)x)\nu(dx)\quad\text{for $B$ 
Borel with $B\not\ni0$},\\
\gamma_t=\int_0^t f(s)\left(\gamma+\int_{\mathbb{R}^d}x\left(\frac{1}{1+|f(s)x|^2}
-\frac{1}{1+|x|^2}\right)\nu(dx)\right)ds.
\end{gather*}
Hence, recalling that $\int_{|x|>1}|x|\nu(dx)<\infty$, we obtain
\begin{align*}
&\int_{|x|>1}|x|\nu_t(dx)=\int_0^t ds\int_{|f(s)x|>1} |f(s)x|\nu(dx)\\
&\qquad\leqslant \int_0^t ds\int_{|x|>1}|f(s)x|\nu(dx)+\int_0^t ds\int_{|x|\leqslant1}
|f(s)x|^2\nu(dx)\\
&\qquad= \int_0^t |f(s)|ds\int_{|x|>1}|x|\nu(dx)+\int_0^t f(s)^2ds
\int_{|x|\leqslant1}|x|^2\nu(dx)\\
&\qquad<\infty.
\end{align*}
Thus $E|Y_t|<\infty$. Now we have
\begin{align*}
EY_t&=\gamma_t+\int_{\mathbb{R}^d}\frac{x\,|x|^2}{1+|x|^2}\nu_t(dx)\\
&=\int_0^t f(s)\left(\gamma+\int_{\mathbb{R}^d}x\left(\frac{1}{1+|f(s)x|^2}
-\frac{1}{1+|x|^2}\right)\nu(dx)\right)ds\\
&\quad+\int_0^t ds\int_{\mathbb{R}^d}\frac{f(s)x\,|f(s)x|^2\nu(dx)}{1+|f(s)x|^2}\\
&=\int_0^t f(s)\left(\gamma+\int_{\mathbb{R}^d}\frac{x\,|x|^2\nu(dx)}{1+|x|^2}\right)ds
=0,
\end{align*}
that is,} $\{Y_t\}$ is a martingale additive process.  \qed

\medskip
{\it Remark on Proposition \ref{mart}}.  If $X^{(\mu)}$ is a martingale 
L\'evy process and if $f\in\mathbf L(X^{(\mu)})$, it is not necessarily true
that $Y_t=\int_0^t f(s)dX_s^{(\mu)}$ is a martingale additive process.
In fact, $E|Y_t|$ may be infinite. 
For example let $X^{(\mu)}$  be a compound Poisson process on $\mathbb{R}^d$
with mean zero.
Then any measurable function $f$ belongs to 
$\mathbf L(X^{(\mu)})$ as Example 4.4 of \cite{S05D} says.
If $\int_0^{t_0} |f(s)|ds=\infty$, then $E|Y_{t_0}|=\infty$, because, choosing
$a>0$ such that $0<\int_{|x|>a} |x|\nu(dx)<\infty$ for the L\'evy 
measure $\nu$ of $X^{(\mu)}$, we have, for the L\'evy 
measure $\nu_{t_0}$ of $Y_{t_0}$, 
{\allowdisplaybreaks
\begin{align*}
\int_{|x|>1}|x|\nu_{t_0} (dx)&=\int_0^{t_0} ds\int_{|f(s)x|>1} |f(s)x|\nu(dx)\\
&\geqslant \int_{|x|>a}|x|\nu(dx)\int_{[0,{t_0}]\cap\{|f(s)|>1/a\}} |f(s)|ds
=\infty,
\end{align*}
which implies that} $E|Y_{t_0}|=\infty$.

\medskip
{\it Remark on martingale additive processes related to Theorem \ref{t5a}}. 
The L\'evy process $X^{(\mu)}$ associated with $\mu$ in Theorems \ref{t5a} is a
martingale, that is, it satisfies $E|X_t^{(\mu)}|<\infty$ and
$EX_t^{(\mu)}=0$. 
Consider the case $d=1$. Let $h(s)$, $D^+$, $D^-$, and $D^0$ be as in the 
proof of Theorem \ref{t5a}.  Thus
\[
h(s)=s^{-1}\gamma +\int_{|x|\geqslant2} s^{-1}x \left(
\frac{1}{1+(s^{-1}x)^2}-\frac{1}{1+x^2}\right) \nu(dx)\quad\text{for }s\geqslant1
\]
and $D^+$, $D^-$, or $D^0$ is the set of $s\geqslant1$ at which $h(s)$ is positive,
negative, or zero, respectively.
Let
{\allowdisplaybreaks
\begin{gather*}
Y_t=\int_0^t s^{-1}1_{[1,\infty)}(s)dX_s^{(\mu)},\\
Y_t^p=\int_0^t s^{-1} 1_{D^p}(s)dX_s^{(\mu)}\quad\text{for }p=+,-,0.
\end{gather*}
Then $\{Y_t\}$, $\{Y_t^+\}$, $\{Y_t^-\}$, $\{Y_t^0\}$, and $\{Y_t^+ +Y_t^-\}$
are martingale additive processes, as is shown in Proposition \ref{mart}.
We can show that
\begin{gather}
Y_t^+\to\infty\quad\text{and}\quad Y_t^-\to-\infty\quad\text{a.s. as }t\to\infty,
\label{Y1}\\
\text{$Y_t$, \,$Y_t^+ +Y_t^-$, and $Y_t^0$ are convergent in $\mathbb{R}$ a.s.\ 
as $t\to\infty$},\label{Y3}\\
E|Y_{\infty-}|=\infty.\label{Y4}
\end{gather}
These are remarkable behaviors. If} $A=0$, then each of $Y_t^+$ and $Y_t^-$ is the
compensated sum of the jumps of $X^{(\mu)}$ in the union of 
some nonrandom time intervals with some nonrandom weights.
For these behaviors it is essential that the L\'evy measure is nonsymmetric
and close to symmetric.  Theorem \ref{t3} (ii) says that martingale compound
Poisson processes with symmetric L\'evy measures do not exhibit this kind of
behaviors.

Proof of \eqref{Y1}--\eqref{Y4} is as follows.
Let the triplet of $Y_t^p$ be $(A_t^p, \nu_t^p, \gamma_t^p)$ for $p=+,-,0$. Then
{\allowdisplaybreaks
\begin{gather*}
A_t^p=\int_1^t s^{-2} 1_{D^p}(s)ds\,A,\\
\nu_t^p(B)=\int_1^t 1_{D^p}(s)ds\int_{\mathbb{R}} 1_B(s^{-1}x)\nu(dx)\quad\text{for $B$ 
Borel set with $B\not\ni0$},\\
\gamma_t^p=\int_1^t  1_{D^p}(s)h(s)ds.
\end{gather*}
Since} $A_t^p$ and $\int_{\mathbb{R}}(x^2\land1)\nu_t^p(dx)$ are bounded and increasing,
$Y_t^p -\gamma_t^p$ is
convergent in probability as $t\to\infty$.  Since it is an additive process,
$Y_t^p -\gamma_t^p$ is convergent a.s.\ also. Since
$\gamma_t^+\to\infty$ and $\gamma_t^-\to -\infty$ (see the proof of Theorem \ref{t5a}),
we obtain \eqref{Y1}.
We have convergence of $Y_t^0$ since $\gamma_t^0=0$. 
Convergence of $Y_t$  comes from the fact that 
$\mu\in\mathfrak D(\Phi_{f_1})$.
Recalling that $Y_t=Y_t^+ +Y_t^- +Y_t^0$, we obtain \eqref{Y3}. 
In order to see \eqref{Y4}, let $\nu_{\infty-}$ be
the L\'evy measure of $Y_{\infty-}$.
Then
{\allowdisplaybreaks
\begin{align*}
&\int_{|x|>1} |x|\nu_{\infty-}(dx)=\int_1^{\infty} ds\int_{|s^{-1}x|>1}
|s^{-1}x|\nu(dx)\\
&\qquad=\int_{|x|>1} |x|\nu(dx)\int_1^{|x|} s^{-1}ds
=\int_{|x|>1} |x|\log |x|\nu(dx),
\end{align*}
which is infinite by virtue of Theorem 2.8 of \cite{S05A} and of the 
fact that} $\mu\not\in \mathfrak D^0(\Phi_{f_1})$.
Hence $E|Y_{\infty-}|=\infty$.  

\section{Applications}

The following results are consequences of Theorems \ref{t1} and \ref{t2}.

\begin{prop}\label{p3.1}
Let $f_1$ and $f_2$ be measurable and $|f_2|\leqslant |f_1|$. 
If $\mathfrak D(\Phi_{f_1})=\mathfrak D^0(\Phi_{f_1})$ or if $\mathfrak D(\Phi_{f_2})=\mathfrak{D}_{\mathrm{e}}
(\Phi_{f_2})$, then $\mathfrak D(\Phi_{f_1})\subset \mathfrak D(\Phi_{f_2})$.
\end{prop}

{\it Proof}. In general we have \eqref{8}. 
Hence it follows from Theorem \ref{t2} that if 
$\mathfrak D(\Phi_{f_1})=\mathfrak D^0(\Phi_{f_1})$, then
\[
\mathfrak D(\Phi_{f_1})=\mathfrak D^0(\Phi_{f_1})\subset \mathfrak D^0(\Phi_{f_2})\subset
\mathfrak D(\Phi_{f_2});
\]
it follows from Theorem \ref{t1} that if $\mathfrak D(\Phi_{f_2})=\mathfrak{D}_{\mathrm{e}}
(\Phi_{f_2})$, then
\[
\mathfrak D(\Phi_{f_1})\subset\mathfrak{D}_{\mathrm{e}}(\Phi_{f_1})\subset\mathfrak{D}_{\mathrm{e}}(\Phi_{f_2})=
\mathfrak D(\Phi_{f_2}),
\]
completing the proof.  \qed

\begin{ex}\label{e3.1}
Let $f_1$ be a locally square-integrable function on $[0,\infty)$
satisfying $f_1(s)\asymp s^{-1/\alpha}$ as $s\to\infty$ with some
$\alpha\in(0,1)\cup(1,2)$.  Let $f_2(s)$ and $f_3(s)$ be measurable and
satisfy $|f_2(s)|\leqslant |f_1(s)|\leqslant|f_3(s)|$.  If $\alpha\in(0,1)$, then 
$\mathfrak D(\Phi_{f_3})\subset\mathfrak D(\Phi_{f_1})\subset \mathfrak D(\Phi_{f_2})$.
If $\alpha\in(1,2)$, then 
$\mathfrak D(\Phi_{f_1})\subset \mathfrak D(\Phi_{f_2})$.

Indeed, we have $\mathfrak D(\Phi_{f_1})=\mathfrak D^0(\Phi_{f_1})=\mathfrak{D}_{\mathrm{e}}(\Phi_{f_1})$
if $\alpha\in(0,1)$, and 
$\mathfrak D(\Phi_{f_1})=\mathfrak D^0(\Phi_{f_1})\subsetneqq\mathfrak{D}_{\mathrm{e}}(\Phi_{f_1})$ if 
$\alpha\in(1,2)$ (Theorem 2.4 of \cite{S05A}).  Hence Proposition \ref{p3.1}
applies.
\end{ex}

\begin{prop}\label{p3.2}
Let $f$ be a locally square-integrable function on $[0,\infty)$
such that there are positive constants $\alpha$, $c_1$, and $c_2$ satisfying
\begin{equation}\label{3.1}
e^{-c_2 s^{\alpha}}\leqslant f(s)\leqslant e^{-c_1 s^{\alpha}}\quad\text{for all large $s$.}
\end{equation}
Then
\begin{equation}\label{3.2}
\begin{split}
&\mathfrak D^0(\Phi_f)=\mathfrak D(\Phi_f)=\mathfrak{D}_{\mathrm{c}}(\Phi_f)=\mathfrak{D}_{\mathrm{e}}(\Phi_f)\\
&\qquad=\left\{\mu\in ID(\mathbb{R}^d)\colon \int_{\mathbb{R}^d}(\log^+ |x|)^{1/\alpha}\mu(dx)<\infty\right\}\\
&\qquad=\left\{\mu\in ID(\mathbb{R}^d)\colon \int_{\mathbb{R}^d}(\log^+ |x|)^{1/\alpha}\nu(dx)<\infty\right\},
\end{split}
\end{equation}
where $\nu$ is the L\'evy measure of $\mu$ and $\log^+ u=(\log u)\lor0$.
\end{prop}

Obviously, in \eqref{3.2}, we can use $(\log(1+|x|))^{1/\alpha}$ for $|x|>1$ in
place of $(\log^+ |x|)^{1/\alpha}$.

\medskip
{\it Proof of Proposition \ref{p3.2}}.  
Let $M=\left\{\mu\in ID(\mathbb{R}^d)\colon \int_{\mathbb{R}^d}(\log^+ |x|)^{1/\alpha}
\mu(dx)<\infty\right\}$. Then $M$ has the last expression in \eqref{3.2}, which
is a consequence of  Theorem 25.3 of \cite{S}.
 Let $f_j(s)=e^{-c_j s^{\alpha}}$, $j=1,2$. Using Theorem 5.15 of \cite{S05D} for
these functions,  we see that
\[
\mathfrak D(\Phi_{f_j})=\mathfrak{D}_{\mathrm{c}}(\Phi_{f_j})=\mathfrak{D}_{\mathrm{e}}(\Phi_{f_j})=M,\qquad j=1,2.
\]
Combined with Proposition \ref{p3} of this paper, the proof of that theorem also
shows that $\mathfrak D^0(\Phi_{f_j})=M$.
Since $f_2(s)\leqslant f(s)\leqslant f_1(s)$ for all large $s$, 
it follows from Theorems \ref{t1} and \ref{t2} that
\[
\mathfrak D^0(\Phi_{f_1})\subset \mathfrak D^0(\Phi_{f})\subset \mathfrak D^0(\Phi_{f_2}),\qquad
 \mathfrak{D}_{\mathrm{e}}(\Phi_{f_1})\subset \mathfrak{D}_{\mathrm{e}}(\Phi_{f})\subset 
\mathfrak{D}_{\mathrm{e}}(\Phi_{f_2}).
\]
Thus $\mathfrak D^0(\Phi_f)= \mathfrak{D}_{\mathrm{e}}(\Phi_f)=M$.  Using \eqref{8},
we also have $\mathfrak D(\Phi_f)= \mathfrak{D}_{\mathrm{c}}(\Phi_f)=M$.  \qed

\medskip
Theorem 5.15 of \cite{S05D} deals with a function $f(s)$ such that 
$f(s)\asymp s^{\beta} e^{-cs^{\alpha}}$, $s\to\infty$, with $\alpha>0$, $\beta\in\mathbb{R}$, and
$c>0$.  This function satisfies \eqref{3.1}.  Thus, if we show Theorem 5.15 of 
\cite{S05D} only for $f(s)=e^{-cs^{\alpha}}$, then the proof of our Proposition 
\ref{p3.2} is obtained and the rest of Theorem 5.15 of 
\cite{S05D} is a consequence of our Proposition \ref{p3.2}.

\begin{ex}\label{e3.2}
Let $f$ be as in Proposition \ref{p3.2}.  If 
$f_2(s)$ and $f_3(s)$ are measurable and
satisfy $|f_2(s)|\leqslant |f(s)|\leqslant|f_3(s)|$, then
$\mathfrak D(\Phi_{f_3})\subset\mathfrak D(\Phi_{f})\subset \mathfrak D(\Phi_{f_2})$.
Use Propositions \ref{p3.1} and \ref{p3.2}.
\end{ex}

\medskip
Let $L_0(\mathbb{R}^d)$ be the class of selfdecomposable distributions on $\mathbb{R}^d$ and let
$L_m(\mathbb{R}^d)$, $m=1,2,\ldots$, be the nested subclasses of $L_0(\mathbb{R}^d)$ studied by
Urbanik \cite{U72,U73} and Sato \cite{S80}. The stochastic integral representation
of $L_0(\mathbb{R}^d)$ given by Wolfe \cite{W82a,W82b}, Jurek and Vervaat \cite{JV83}, 
and others is in the form $\Phi_f$ with $f(s)=e^{-s}$.  Further, the representation 
of $L_m(\mathbb{R}^d)$ for $m=1,2,\ldots$ 
given by Jurek \cite{J83} can be rewritten in the form $\Phi_f$
with $f(s)=e^{-cs^{1/(m+1)}}$.  Hence we can apply Proposition \ref{p3.2}
to those cases.
Further applications related to \cite{BMS04a} are in progress.


\bigskip
\begin{thebibliography}{22}

\bibitem{BMS04a} Barndorff-Nielsen, O.E., Maejima, M., 
Sato, K. (2006) Some classes
of infinitely divisible distributions admitting 
stochastic integral representations,
{\it Bernoulli\/} {\bf12}, 1--33.

\bibitem{CS05} Cherny, A, Shiryaev, A. (2005) On stochastic integrals
up to infinity and predictable criteria for integrability, {\it Lecture
Notes in Math.}, Springer, {\bf1857}, 165--185.

\bibitem{J83}  Jurek, Z.J. (1983) The class $L_{m}(Q)$ of
probability measures on Banach spaces, {\it Bull.\ Polish 
Acad.\ Sci.\ Math.\/} {\bf31}, 51--62.

\bibitem{JV83} Jurek, Z.J., Vervaat, W. (1983) An integral representation
for selfdecomposable Banach space valued random variables, {\it Zeit.\ 
Wahrsch.\ Verw.\ Gebiete} {\bf62}, 247--262.

\bibitem{KW92}Kwapie\'n, S., Woyczy\'nski, W.A. (1992) 
{\it Random Series and 
Stochastic Integrals: Single and Multiple}, Birkh\"auser, 
Boston.

\bibitem{RR89} Rajput, B., Rosinski, J. (1989) Spectral 
representations of infinitely 
divisible processes, {\it Probab.\ Theory Related Fields\/} 
{\bf82}, 451--487.

\bibitem{S80}Sato, K. (1980) Class $L$ of multivariate distributions 
and its subclasses,
{\it J.\ Multivar.\ Anal.\/} {\bf10}, 207--232.

\bibitem{S}
 Sato, K. (1999) {\it L\'evy Processes and Infinitely 
Divisible Distributions}, 
Cambridge Univ.\ Press, Cambridge.

\bibitem{S04}
 Sato, K. (2004) Stochastic integrals in additive processes 
and application to
semi-L\'evy processes, {\it Osaka J. Math.\/} {\bf41}, 211--236.

\bibitem{S05D}
 Sato, K. (2005) Additive processes and  stochastic integrals, preprint.

\bibitem{S05A}
 Sato, K. (2006) Two families of improper stochastic integrals with respect 
to L\'evy processes, {\it Alea, Latin American 
Journal of Probability and
 Mathematical Statistics} {\bf1}, 47--87. http://alea.impa.br/english/

\bibitem{U72}  Urbanik, K. (1972) Slowly varying sequences of random variables,
{\it Bull.\ Acad.\ Polon.\ Sci.\ S\'er.\ Sci.\ Math.\ Astronom.\ Phys.} {\bf20}, 
678--682.

\bibitem{U73}  Urbanik, K. (1973) Limit laws for sequences of normed sums
satisfying some stability conditions, {\it Multivariate Analysis--III}
(ed. Krishnaiah, P.\thinspace R., Academic Press, New York), 225--237.

\bibitem{UW67} Urbanik, K., Woyczy\'nski, W.A. (1967) 
Random integrals and Orlicz spaces,
{\it Bull.\ Acad.\ Polon.\ Sci.} {\bf15}, 161--169.

\bibitem{W82a} Wolfe, S.J. (1982) A characterization of certain stochastic 
integrals, (Tenth Conference on Stochastic Processes and Their Applications, 
Contributed Papers) {\it Stoch.\ Proc.\ Appl.} {\bf12}, 136--136.

\bibitem{W82b} Wolfe, S.J. (1982) On a continuous analogue of the stochastic difference
equation $X_n=\rho X_{n-1}+B_n$, {\it Stoch.\ Proc.\ Appl.} {\bf12}, 301--312.
\end{thebibliography}
\end{document}